\documentclass[11pt, twocolumn]{article}
\pdfoutput=1
\usepackage{latexsym}
\usepackage{amsfonts}
\usepackage{amssymb,amsmath}
\usepackage{graphicx}
\begin{document}
\date{}
\title{
\vspace{-1cm}
{\LARGE
\textbf{A Spectral Method for Solving the Cauchy Problem}}
\footnote{This is a corrected version of the paper 
\em
``A spectral method for Cauchy problem solving''
\em published in Proceedings of MCME International Conference, Problems in Modern 
Applied Mathematics, N. Mastorakis, ed., WSES Press, 2000, pp.~227-232.
Alterations were made in the explanatory part of the paper, whereas core mathematical constructions were copied exactly. Footnotes and an Appendix were added when correcting.}
}
\author{
V.S. CHELYSHKOV\\
\\
Institute of Hydromechanics, NAS,\\
8/4 Zhelyabov St., Kiev, 03057\\
UKRAINE\\
\\
\\
\normalsize
\hspace{-0.3cm}
\em Abstract:\em 
\hspace{2pt}A new approach for integration of the initial value problem for ordinary differential equations is\\
\normalsize
\hspace{-0.3cm}
suggested.
The algorithm is based on approximation of 
the solution by a system of functions
that contains\\
\hspace{-12cm}
\normalsize orthogonal exponential polynomials.\\
\vspace{10pt}
\\
\normalsize 
\hspace{-0.4cm}
\normalsize 
\hspace{-16.1cm}
integration. 
}
\maketitle
\vspace{-0.5cm}
\textbf{\Large \hspace{-0.6cm} 1 Introduction}

\hspace{-0.55cm}
The theory of integration of the initial value problem for ordinary 
differential equations is well developed. It includes various and effective methods (\cite{Gear}, \cite{Stoer}, for example) that are widely used \cite{Press}. 
A new approach for solving the Cauchy problem is introduced in this paper. It is based on approximation of the solution by a sequence of functions and deriving $n$ nonlinear functional equations from the initial value problem.  
The procedure does not result in consecutive elimination of low order terms of the local error when increasing the order of spectral approximation. Thus, presented discretization  is different from known implicit Runge-Kutta methods.  

We employ the nonlinear approximation for constructing a one-step explicit algorithm. Description of the algorithm is illustrated by an example on numerical integration of a large-scale nonlinear dynamical system.
\\
\\
\\
\textbf{\Large 2 Discretization of the Problem}

\hspace{-0.5cm}
Let vector $\mbox{\boldmath ${\it Y}$}_{I}$ and vector functions \mbox{\boldmath ${\it Y}$} and \mbox{\boldmath ${\it F}$} have dimension $N$. We consider the problem
\[
\mbox{\boldmath ${\it Y}$}^\prime(T)=\mbox{\boldmath ${\it F}$}
(T,\mbox{\boldmath ${\it Y}$}(T)),
\]
\begin{equation}
\label{eq1}
\mbox{\boldmath $\it Y$}(T_I)=
\mbox{\boldmath $\it Y$}_{I},\quad
T\in\left[T_I,T_E\right],\quad T_I < T_E,
\end{equation}
supposing that the function \mbox{\boldmath ${\it F}$} has properties that are necessary for the following constructions.

Let $n\in {\mathbb N}$, $h=T_E-T_I$, $\alpha_n=h/\lambda_{nn}$, and $\lambda_{nn}$ is a parameter that will be defined below. Making use of the substitutions
\begin{equation}
\label{eq2}
t=(T-T_I)/\alpha_n,\quad
\mbox{\boldmath ${\it y}$}(t)=\mbox{\boldmath ${\it Y}$}
(T_I+\alpha_{n}t),
\end{equation}
\begin{equation}
\label{eq3}
\mbox{\boldmath ${\it f}$}(t,\mbox{\boldmath ${\it y}$})=
\alpha_n\mbox{\boldmath ${\it F}$}(T_I+\alpha_{n}t,
\mbox{\boldmath ${\it Y}$}(T_I+\alpha_{n}t))
\end{equation}
we reduce the original problem to the following
\begin{equation}
\label{eq4}
\mbox{\boldmath ${\it y}$}^{\prime}(t)=
\mbox{\boldmath ${\it f}$}(t,\mbox{\boldmath ${\it y}$}(t)),\quad
\mbox{\boldmath ${\it y}$}(0)=\mbox{\boldmath ${\it Y}$}_I,\quad
t\in\left[0,\lambda_{nn}\right].
\end{equation}
We look for an approximation 
$\mbox{\boldmath ${\it y}$}_{n}(t)\approx\mbox{\boldmath ${\it y}$}(t)$ to problem (\ref{eq4}) in 
the form
\begin{equation}
\label{eq5}
\mbox{\boldmath ${\it y}$}_{n}(t)=
\mbox{\boldmath ${\it Y}$}_I+
\mbox{\boldmath ${\it a}$}_{n0}t+
\sum\limits_{j=1}^{n}
\mbox{\boldmath ${\it a}$}_{nj}
\mbox{${\cal S}$}_{nj}(\beta_n, t)
\end{equation}
\begin{equation}
\label{eq6}
\mbox{${\cal S}$}_{nj}(\beta_n, t)=\int_0^t
\mbox{${\cal E}$}_{nj}(\beta_n t)\mbox{d}t
\end{equation}
that satisfies the initial conditions in (\ref{eq4}). The unknown vector coefficients  
$\mbox{\boldmath ${\it a}$}_{n0}$, 
$\mbox{\boldmath ${\it a}$}_{nj}$ 
have dimension
$N$, $\beta_n=1$,  $\mbox{${\cal E}$}_{nj}(t)$ are 
the polynomials of exponents, which are orthogonal on the semi-axis with the weight function 1. The polynomials can be defined in the following 
recurrent way \cite{VC1}:
\[
{\cal E}_{nn}(t)=\mbox{exp}(-nt),
\]
\[
{\cal E}_{n,n-1}(t)=
(2n-1) \mbox{exp}(-(n-1)t)-2n \mbox{exp}(-nt),
\]
\[
a_{nj}{\cal E}_{n,j-1}(t) =(b_{nj}
\mbox{exp}(t) - c_{nj} ){\cal E}_{nj}(t) -d_{nj}{\cal E}_{n,j+1}(t),
\]
\[
a_{nj} =(2j+1)(n+j)(n-j+1),
\]
\[
b_{nj} =(2j-1)2j(2j+1),
\]
\[
c_{nj} =4j(n^2+j^2+n),
\]
\[
d_{nj} =(2j-1)(n-j)(n+j+1),
\]
$j=n-1,n-2,\ldots, 2$. The recurrence relations also can be extended to $j=1$, and this results in definition of the polynomial ${\cal E}_{n0}(t)$ which has the 
zeros $\lambda_{ns}$, $s=1,\ldots, n$, $\lambda_{nn}$ is the maximum non-trivial zero. The polynomial ${\cal E}_{n0}(t)$ is orthogonal to the system of functions
$\mbox{\boldmath ${\cal E}$}_{n}(t)=\{{\cal E}_{nj}(t)\}_{j=1}^n$,
but it is not integrable on the semi-axis. 

Obviously, $\mbox{lim}_{t\rightarrow\infty}\mbox{\boldmath ${\cal E}$}_{n}(t)=0$, and one can  come to the conclusion that $\mbox{\boldmath ${\cal E}$}_{n}(t)$ is complete in $L_{2}[0,\infty)$.  
Hence, nice properties of the system on the semi-axis may serve a useful purpose in approximation of a continuous function on a finite interval if the system is completed by unity.   
Accordingly, the sequence of functions 
$1,t,\mbox{\boldmath ${\cal E}$}_{n}(t)$ may be employed for 
approximation of a continuously differentiable function on an interval. Such a deduction leads to approximation (\ref{eq5}), (\ref{eq6}), and the integral in the right-hand side of (\ref{eq6}) is calculated by the formula
\[
\int_0^t
{\cal E}_{nj}(\beta_n t)\mbox{d}t=
\]
\[\frac{1}{\beta_n j}
\left(1-{\cal E}_{nj}(\beta_n t)-2\sum\limits_{l=j+1}^{n}{\cal E}_{nl}(\beta_n t)\right).
\]
Satisfying the initial condition for the derivative of the 
approximation $\mbox{\boldmath ${\it y}$}_{n}(t)$ in (\ref{eq5}) and making use of the polynomial values ${\cal E}_{nj}(0)=(-1)^{n-j}$ one can immediately show that
\begin{equation}
\label{eq7}
\mbox{\boldmath ${\it y}$}_{n}(t)=
\mbox{\boldmath ${\it Y}$}_I+
\mbox{\boldmath ${\it f}$}_{0}t+ 
\sum\limits_{j=1}^{n} 
\mbox{\boldmath ${\it a}$}_{nj}
({\cal S}_{nj}(\beta_n, t)-(-1)^{n-j}t),
\end{equation}
where
$\mbox{\boldmath ${\it f}$}_{0}=
\mbox{\boldmath ${\it f}$}(0,
\mbox{\boldmath ${\it Y}$}_{I})$.
 
To find the coefficients
$\mbox{\boldmath ${\it a}$}_{nj}$ 
in (\ref{eq7}) we 
consider the discrete form 
\begin{equation}
\label{eq8}
\mbox{\boldmath ${\it y}$}_{n}^{\prime}(t_{ns})=\mbox{\boldmath ${\it f}$}
(t_{ns},\mbox{\boldmath ${\it y}$}_{n}(t_{ns}))
\end{equation}
of the ordinary differential
equation in (\ref{eq4}) at the collocation points $t_{ns}=\lambda_{ns}/\beta_n$, $s=1,\ldots,n$.
We are going to reduce equations (\ref{eq8}) to a system of functional equations, and we quote below a few more properties of polynomials ${\cal E}_{nj}(t)$ that are necessary for computations.

One can observe that the system of functions
$(\mbox{\boldmath ${\cal E}$}_{n}\cup{\cal E}_{n0})(t)$ generates the Gauss-type quadratures for exponents on
the semi-axis.
This results in the discrete form of orthogonality of the polynomials
\footnote[1]{{\sc V.~S. Chelyshkov}, A variant of spectral method in the theory of hydrodynamic stability, 
{\em Hydromechanics} ({\em Gidromekhanika}), N 68, 1994, pp.~105--109, (in Russian).}
:
\[
\sum\limits_{s=1}^n \rho_{ns}{\cal E}_{nj}(\lambda_{ns}){\cal E}_{nl}(\lambda_{ns})=
\delta_{jl}/(j+l),
\]
\begin{equation}
\label{eq9}
\rho_{ns}=1/\left(2\sum\limits_{m=1}^{n}m{\cal E}_{nm}^2(\lambda_{ns})\right),
\end{equation}
$j,l=1,2,\ldots,n$, $\delta_{jl}$ is the Kronecker delta, $\lambda_{ns}$ and $\rho_{ns}$ are the abscissas and weights of the quadrature. 
Also, the formula 
\[
\sum\limits_{s=1}^n \rho_{ns}{\cal E}_{nl}(\lambda_{ns})=1/l
\]
takes place. It stands for the integral of ${\cal E}_{nl}(t)$ on the semi-axis.
It may be noticed that the formula for the weights in (\ref{eq9}) is similar to that one in \cite{Lanczos} for the original Gauss quadrature.  The Gauss rule for exponents also follows straightforward from the substitutions
\begin{equation}
\label{eq10}
z_{ns}=1-2\mbox{exp}(-\lambda_{ns}),\quad
w_{ns}=2\rho_{ns}\mbox{exp}(-\lambda_{ns}),
\end{equation}     
were $z_{ns}$ and $w_{ns}$ are the abscissas and weights of the original Gauss quadrature on the interval $[-1,1]$. 
 
Substituting approximation (\ref{eq7}) to equations (\ref{eq8}), multiplying the equations by 
$\rho_{ns}{\cal E}_{nl}(\lambda_{ns})$, adding them, making use of the properties
described, and 
inverting the matrix in the left-hand side of the equality developed one can obtain the following equations 
\begin{equation}
\label{eq11}
\mbox{\boldmath ${\it a}$}_{nj}=2\sum\limits_{l=1}^n A_{jl}l\sum\limits_{s=1}^n\rho_{ns}{\cal E}_{nl}(\lambda_{ns})\mbox{\boldmath ${\it f}$}_{ns}
-2(-1)^{n}\mbox{\boldmath ${\it f}$}_0,
\end{equation}
where 
\[
\mbox{\boldmath ${\it f}$}_{ns}\equiv
\mbox{\boldmath ${\it f}$}(t_{ns},\mbox{\boldmath ${\it y}$}_{n}(t_{ns})),
\]
\[
A_{jl}=
\left \{
\begin{array}{c}
(-1)^{l}2, \\ [1ex]
\,\,\,\,\,\,\,\,\,\,-1,   \\ [1ex]
\,\,\,\,\,\,\,\,\,\,\,\,\,\,3,
\end{array}
\right.
\begin{array}{c}
j\neq l, \\ [1ex]
j=l,	\\ [1ex]
j=l,
\end{array}
\begin{array}{c}
	\\ [1ex]
l  \mbox{  odd,	}\\ [1ex]
l  \mbox{  even.}
\end{array}
\]
Equation (\ref{eq11}) together with approximation (\ref{eq7})  at $t=~t_{ns}$ form the system of $N\times n$ functional 
equations for the components of the vectors $\mbox{\boldmath ${\it a}$}_{nj}$.
However, we will not examine such a spectral form of the discretization, but return to the original variables. 

Substituting (\ref{eq11}) to (\ref{eq7}) we represent the approximation in the form
\begin{equation}
\label{eq12}
\mbox{\boldmath ${\it y}$}_{n}(t)=
\mbox{\boldmath ${\it Y}$}_{\hspace{-1mm}I}+
\mbox{\boldmath ${\it f}$}_{0}
{\cal Q}_{n0}(t)
+ 
\sum\limits_{s=1}^{n}
\mbox{\boldmath ${\it f}$}_{ns} 
{\cal Q}_{ns}(t)
,
\end{equation}
where the functions 
in the right hand side are 
\begin{equation}
\label{eq13}
{\cal Q}_{n0}(t)=(-1)^n\left( t-2\sum\limits_{j=1}^n {\cal S}_{nj}(\beta_n,t)\right),
\end{equation}
\[
\hspace{-23mm}
{\cal Q}_{ns}(t)=2\rho_{ns}\sum\limits_{l=1}^nl{\cal E}_{nl}(\lambda_{ns})
\times
\]
\begin{equation}
\label{eq14}
\hspace{24mm}\left(\sum\limits_{j=1}^nA_{jl}{\cal S}_{nj}(\beta_n,t)-(-1)^lt\right).
\end{equation}
Approximation (\ref{eq12}) -- (\ref{eq14}) gives the opportunity to evaluate $\mbox{\boldmath ${\it y}$}_{n}(t)$ for any 
current $t$ if the values of $\mbox{\boldmath ${\it y}$}_{n}(t_{ns})$ are known.

Equating (\ref{eq12}) at the points $t=t_{np}$, $p=1, \ldots, n$ 
and making use of substitutions (\ref{eq2}), (\ref{eq3}),
we finally obtain the discrete analogue of initial value problem (\ref{eq1}) in the form 
\[
T_{np}=T_I+\nu_{np}h,\quad 
\mbox{\boldmath ${\it Y}$}_{\hspace{-1mm}np}=\mbox{\boldmath ${\it Y}$}(T_{np}),
\]
\begin{equation}
\label{eq15}
\mbox{\boldmath ${\it Y}$}_{\hspace{-1mm}np}=\mbox{\boldmath ${\it Y}$}_{\hspace{-1mm}I}+ \sigma_{np0}h\mbox{\boldmath ${\it F}$}(T_I, \mbox{\boldmath ${\it Y}$}_{\hspace{-1mm}I})+h\sum\limits_{s=1}^n\sigma_{nps}\mbox{\boldmath ${\it F}$}(T_{ns}, \mbox{\boldmath ${\it Y}$}_{\hspace{-1mm}ns}),
\end{equation}
where the coefficients in (\ref{eq15}) are calculated as follows 
\[
\nu_{np}=t_{np}/\lambda_{nn},
\]
\[
\sigma_{np0}={\cal Q}_{n0}(t_{np})/\lambda_{nn},
\quad
\sigma_{nps}={\cal Q}_{ns}(t_{np})/\lambda_{nn}.
\]

Although discretization (\ref{eq15}) resembles an implicit Runge-Kutta method, the procedure developed is not subjected to the basic idea of the method. In fact, low order terms in the expansion of the approximate solution in the Taylor series in $h$ are not eliminated when increasing $n$. 

If $\mbox{\boldmath ${\it F}$}(T,\mbox{\boldmath ${\it Y}$}(T))$ is a linear function of $\mbox{\boldmath ${\it Y}$}$ then the precise solution of problem (\ref{eq15}) may be found
\footnote[2]{If $\mbox{\boldmath ${\it F}$}(T,\mbox{\boldmath ${\it Y}$}(T))\equiv\mbox{\boldmath ${\it F}$}(T)$, $T_I=0$, and $h=1$ then equality (\ref{eq15}) contains the quadrature rule that is exact for the integrands $\{1, \mbox{exp}(-j\lambda_{nn}T)\}$, $j=1, \ldots, n$ on the interval $[0,1]$. 
}
; in the general case, it is a difficult problem that requires application of methods for solving nonlinear functional equations. 
We will now describe an explicit algorithm that provides approximation to the solution of the initial value problem.
\\
\\
\\
\textbf{\Large 3 A Recurrence Algorithm}

\hspace{-0.5cm}
It follows from the first formula in (\ref{eq10}) that the zeros $\lambda_{ns}$ of the polynomial ${\cal E}_{n0}(t)$ cluster near $t=0$, when $n$ increases. This indicates that computation can originate in the lowest order of approximation and then can be extended from point to point by drawing in the polynomials of successively increasing degree. 

We introduce the recurrence index, $k$, $k=1, \ldots, n$ and substitute $k$ for $n$ in the previous constructions, where it is necessary. We consider successively elongated subintervals $[0, \lambda_{nk}]$ and  collocation points $t_{ks}\equiv t_{ks}(n)$, $s=1, \ldots, k$ on them such that $t_{kk}=\lambda_{nk}$. We meet the requirements of the discrete orthogonality for the points $t_{ks}$ on each of the subintervals in two steps. At first, we consider that $\beta_k \equiv \beta_k(n)$ for $k<n$ and determine $\beta_k$ by satisfying the condition $\beta_k t_{kk}=\lambda_{kk}$, so $\beta_k=\lambda_{kk}/\lambda_{nk}$. Next, we subject the
choice of $t_{ks}$ for $s<k$ to the condition $\beta_k t_{ks}=\lambda_{ks}$, and therefore $t_{ks}=\lambda_{ks}/\beta_{k}$.

We suppose that the integration step of size $h$ is sufficiently 
small and reduce approximation (\ref{eq12}) to the explicit form
\begin{equation}
\label{eq16}
\mbox{\boldmath ${\it y}$}_{k}(t)=
\mbox{\boldmath ${\it Y}$}_{\hspace{-1mm}I}+
\mbox{\boldmath ${\it f}$}_{0}
{\cal Q}_{k0}(t)
+ 
\sum\limits_{s=1}^{k}
\mbox{\boldmath ${\it g}$}_{ks} 
{\cal Q}_{ks}(t)
,
\end{equation}
\begin{equation}
\label{eq17}
\mbox{\boldmath ${\it g}$}_{ks}=
\left \{
\begin{array}{c}
\mbox{\boldmath ${\it y}$}^{\prime}_{k-1}(t_{ks}), \\ [1ex]
\mbox{\boldmath ${\it f}$}(t_{ks},\mbox{\boldmath ${\it y}$}_{k-1}(t_{ks})),   \\ [1ex]
\end{array}
\right.
\begin{array}{c}
\hspace{-3mm}s=1, \ldots, k-1, \\ [1ex]
s=k.	\\ [1ex]
\end{array}
\end{equation}
Thus, we use extrapolation outside of the interval $[0,t_{k-1,k-1}]$ for $s=k$. 

A function, $\mbox{\boldmath ${\it y}$}_{0}(t)$, is required to initiate calculating sequence (\ref{eq16}), (\ref{eq17}). The first two terms in approximation (\ref{eq7}) suggest themselves, and we define 
\[
\mbox{\boldmath ${\it y}$}_{0}(t)=
\mbox{\boldmath ${\it Y}$}_{\hspace{-1mm}I}+
\mbox{\boldmath ${\it f}$}_{0}{\cal R}_{00}(t),\quad
{\cal R}_{00}(t)=t. 
\]
For $k=1$ we have
\[
\mbox{\boldmath ${\it y}$}_{0}(t_{11})=
\mbox{\boldmath ${\it Y}$}_{\hspace{-1mm}I}+
{\cal R}_{00}(t_{11})\mbox{\boldmath ${\it f}$}_{0},
\]
\[
\mbox{\boldmath ${\it g}$}_{11}=
\mbox{\boldmath ${\it f}$}(t_{11},
\mbox{\boldmath ${\it y}$}_{0}(t_{11})),
\]

\[
\mbox{\boldmath ${\it y}$}_{1}(t)=
\mbox{\boldmath ${\it Y}$}_{\hspace{-1mm}I}+
\mbox{\boldmath ${\it f}$}_{0}{\cal R}_{10}(t)+
\mbox{\boldmath ${\it g}$}_{11}{\cal R}_{11}(t),
\]
where
\[
{\cal R}_{10}(t)={\cal Q}_{10}(t),
\quad
{\cal R}_{11}(t)={\cal Q}_{11}(t).
\]
If $n=1$, then $\mbox{\boldmath ${\it y}$}_{1}(t_{11})$ also should be calculated.

Making use of (\ref{eq13}), (\ref{eq14}) and (\ref{eq16}) we calculate the derivatives in (\ref{eq17}), and for $1<k\leq n$ the recurrence relations are as follows
\[
\gamma_{ksr}={\cal R}^{\prime}_{k-1,r}(t_{ks}),
\]
\[
s=1, \ldots, k-1,
\quad
r=0, \ldots, k-1;
\]
\[
\mbox{\boldmath ${\it y}$}^{\prime}_{k-1}(t_{ks})=
\gamma_{ks0}\mbox{\boldmath ${\it f}$}_{0}+
\sum\limits_{r=1}^{k-1}\gamma_{ksr}
\mbox{\boldmath ${\it g}$}_{rr},
\quad
s=1, \ldots, k-1,
\]
\[
\mbox{\boldmath ${\it y}$}_{k-1}(t_{ks})=
\mbox{\boldmath ${\it Y}$}_{\hspace{-1mm}I}+
{\cal R}_{k-1,0}(t_{kk})
\mbox{\boldmath ${\it f}$}_{0}+
\sum\limits_{r=1}^{k-1}{\cal R}_{k-1,r}(t_{kk})
\mbox{\boldmath ${\it g}$}_{rr},
\]
\[
\mbox{\boldmath ${\it g}$}_{kk}=
\mbox{\boldmath ${\it f}$}(t_{kk},
\mbox{\boldmath ${\it y}$}_{k-1}(t_{kk})),
\]
\[
\mbox{\boldmath ${\it y}$}_{k}(t)=
\mbox{\boldmath ${\it Y}$}_{\hspace{-1mm}I}+
\mbox{\boldmath ${\it f}$}_{0}{\cal R}_{k0}(t)+
\sum\limits_{r=1}^{k}
\mbox{\boldmath ${\it g}$}_{rr}{\cal R}_{kr}(t),
\]
where
\[
{\cal R}_{k0}(t)={\cal Q}_{k0}(t)+{\cal G}_{k0}(t)
\]
\[
{\cal R}_{kr}(t)={\cal G}_{kr}(t),
\quad
r=1, \ldots, k-1;
\]
\[
{\cal R}_{kk}(t)={\cal Q}_{kk}(t),
\]
and
\[
{\cal G}_{kr}(t)=
\sum\limits_{s=1}^{k-1}\gamma_{ksr}{\cal Q}_{ks}(t),
\quad
r=0, \ldots, k-1.
\]
If $k=n$, then $\mbox{\boldmath ${\it y}$}_{n}(t_{nn})$
also should be calculated.

Making use of substitutions (\ref{eq2}), (\ref{eq3}) and of previous computations we finally obtain the discretization in the original variables
\begin{equation}
\label{eq18}
\mbox{\boldmath ${\it K}$}_{\hspace{-1mm}n0}=
h\mbox{\boldmath ${\it F}$}(T_{I},
\mbox{\boldmath ${\it Y}$}_{\hspace{-1mm}I}),
\end{equation}
\begin{equation}
\label{eq19}
\mbox{\boldmath ${\it K}$}_{\hspace{-1mm}np}=
h\mbox{\boldmath ${\it F}$}\left(T_{I}+\nu_{np}h, 
\mbox{\boldmath ${\it Y}$}_{\hspace{-1mm}I}+
\sum\limits_{s=0}^{p-1}\mu_{nps}
\mbox{\boldmath ${\it K}$}_{\hspace{-1mm}ns}\right),
\end{equation}
\begin{equation}
\label{eq20}
\mbox{\boldmath ${\it Y}$}_{\hspace{-1mm}n}(T_I+h)=
\mbox{\boldmath ${\it Y}$}_{\hspace{-1mm}I}+ 
\sum\limits_{s=0}^{n}\sigma_{nns}
\mbox{\boldmath ${\it K}$}_{\hspace{-1mm}ns},
\end{equation}
\[
\mu_{nps}={\cal R}_{p-1,s}(t_{pp})/\lambda_{nn},\quad
p=1, \ldots, n.
\]
The algorithm developed requires $n+1$ calculations of the right-hand side 
of problem (\ref{eq1}).
Formulas (\ref{eq18})~ -- ~(\ref{eq20}) are the Runge formulas, but the algorithm does not represent an explicit Runge-Kutta method. Rather, it can be characterized as the explicit Euler method which residual is reduced by the spectral component of the approximation.

Algorithm (\ref{eq18}) -- (\ref{eq20}) may have favorable properties of stability and monotonicity at the expense of precision.  
\\
\\
\\
\textbf{\Large 4 Examples}

\hspace{-0.5cm}
First, we examine the simplest explicit algorithm. On putting $n=1$ in (\ref{eq18}) -- (\ref{eq20}) we find
\[
\mbox{\boldmath ${\it K}$}_{\hspace{-1mm}10}=
h\mbox{\boldmath ${\it F}$}(T_{I},
\mbox{\boldmath ${\it Y}$}_{\hspace{-1mm}I}),
\]
\[
\mbox{\boldmath ${\it K}$}_{\hspace{-1mm}11}=
h\mbox{\boldmath ${\it F}$}\left(T_{I}+h, 
\mbox{\boldmath ${\it Y}$}_{\hspace{-1mm}I}+
\mbox{\boldmath ${\it K}$}_{\hspace{-1mm}10}\right),
\]
\[
\mbox{\boldmath ${\it Y}$}_{\hspace{-1mm}1}(T_I+h)=
\mbox{\boldmath ${\it Y}$}_{\hspace{-1mm}I}+ 
\left(-1+\frac{1}{\mbox{ln}2}\right)
\mbox{\boldmath ${\it K}$}_{\hspace{-1mm}10}+
\left(2-\frac{1}{\mbox{ln}2}\right)
\mbox{\boldmath ${\it K}$}_{\hspace{-1mm}11}.
\]
The algorithm has the same order as the Euler method, but requires two calculations of the right-hand side of the problem. 
The local error at $T=T_{I}+h$ is
\begin{equation}
\label{eq21}
\mbox{\boldmath ${\it Y}$}-
\mbox{\boldmath ${\it Y}$}_{\hspace{-1mm}1}\approx-0.0573
\left(\frac{\partial{\mbox{\boldmath ${\it F}$}}}{\partial t}+
\mbox{\boldmath ${\it F}$}
\frac{\partial{\mbox{\boldmath ${\it F}$}}}
{\partial \mbox{\boldmath ${\it Y}$}}
\right)_{T=T_I}h^2+O(h^3).
\end{equation}
For the Euler method, the constant in (\ref{eq21}) is equal to 0.5. 
For the differential equation $y^{\prime}=-\gamma y$ with $\gamma >~0$, the spectral algorithm holds monotonicity if $h<1.7943/\gamma$, while for the Euler method $h<1/\gamma$.

We elaborated codes for algorithm (\ref{eq18}) -- (\ref{eq20}). Calculation in successively elongated subintervals involves exponential polynomials of increasing degree, and approximations of degree $n-1$ and $n$ can be compared at the end of a step. The comparison was used for fulfillment certain conditions and, as the result, for implementation of adaptive step-size control.        

We examined some algorithms of degree $n\leq16$, and  the next example represents a test on numerical integration of a large-scale nonlinear dynamical system. The system was extracted from an initial~-~ boundary value problem for the Navier-Stokes equations that describes evolution of two dimensional disturbances in the laminar boundary layer near a flat plate \cite{VC3}, \cite{VC2}. Calculations were performed for the Reynolds number $\mbox{R}=10^6$. The dynamical system has $N=1220$ degrees of freedom, and the initial value problem was integrated over a long time interval that displays complete development of disturbances. An example of simulation of two dimensional statistically steady flow for $n=16$ is shown in Figure~ 1 on a representative time interval. Similar graphs also were obtained by a Runge-Kutta method.
\begin{figure}[htbp]
\centerline{\includegraphics{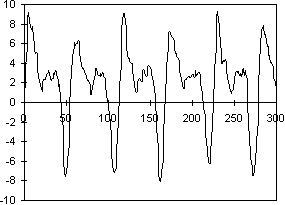}}
\end{figure}

\vspace{-7pt}
\hspace{22pt}
\footnotesize{Figure 1.  Dependence of skin friction on time}
\normalsize

\vspace{12pt}

Since the explicit algorithm is the result of approximation of the solution by a sequence of functions, the solution  can be interpolated in the integration interval and extrapolated outside of it. This could lead to constructing algorithms with translation. The collocation points are distributed non-uniformly on the interval, and the maximum distance between them is proportional to $\lambda_{nn} - \lambda_{n-1,n}$. Thus, an algorithm may include the shift by $h_{\lambda}=(1-\lambda_{n-1,n}/\lambda_{nn})h$, interpolation and extrapolation of the solution at the shifted points, and calculation of the right hand side of the problem at the end point of partially translated interval. 
\\
\textbf{\Large 5 Conclusion}

\hspace{-0.5cm}
A spectral approach for integration of the initial value problem is developed in this paper. The approach results in an implicit procedure and explicit algorithm. The algorithm was tested, and it is comparable with classical explicit methods. In addition, the explicit algorithm may provide a trial value for solving stiff problems.
\\
\\
\\
\textbf{\Large 6 Acknowledgment}

\hspace{-0.5cm}
Author expresses his gratitude to Prof. R. Greechie for a discussion which improved the paper. This work was supported by STCU grant No 473.
\\
\\
\\
\textbf{\Large 7 Appendix}
\footnote[3]{ The author is presently an employee of Lane College, TN, USA.}

\footnotesize{
\hspace{-0.5cm}
Let $N=1$. We consider the stability function (\footnote[4]{{\sc E.~Hairer and G.~Wanner},
{\em Solving Ordinary Differential Equations II}, Springer, Berlin, 1991. },~page 41),
$R_n(z)$, of explicit algorithm (\ref{eq18}) -- (\ref{eq20})  in the complex plane $z$. For $n=1$ and $n=2$ we find
\[
R_{1}(z)\approx1+z+0.557305z^2,
\]
\[
R_{2}(z)\approx1+z+0.533954z^2+0.098846z^3.
\]s
This indicates that convergence of $R_{n}(z)$ to $e^z$ is fairly slow, and the stability domain expands along the negative part of the axis $\Im(z)=0$.

Let $n=1$. For implicit procedure (\ref{eq15}) we find 
\[
R_{1}(z)=\frac{1+Az}{1-(1-A)z},
\]
$A=1/\mbox{ln}2-1\approx 0.442695$, the procedure is A-stable,
and $\lim_{z\rightarrow -\infty}R_{1}(z)\approx-0.794349$.

Alas, for $n=2$ procedure (\ref{eq15}) is not A-stable. Thus, approximation (\ref{eq12}) - (\ref{eq14}) primarily serves as the basis for constructing the explicit algorithm, and discretization (\ref{eq15}) cannot be recommended for solving stiff problems if $n>1$. 

Below we outline two approaches that lead, respectively,  to A- and L-stable procedures.

First, we employ  the same type of approximation, but substitute special sequences of orthogonal exponential polynomials to it.

Let polynomials of exponents ${\cal E}_{nk}^{(\alpha, \beta)}(t)$ be orthogonal sequences on the semi-axis 
with respect to the weight function $(\exp(-t))^{\alpha}(1-\exp(-t))^{\beta}$
such that
\[
{\cal E}_{nk}^{(\alpha, \beta)}(-\mbox{ln}x)=\frac{x^{-\alpha-k}(1-x)^{-\beta}}{(n-k)!}\frac{\mbox{d}^{n-k}}{\mbox{d}x^{n-k}}
(x^{\alpha+n+k}(1-x)^{\beta+n-k}).
\]
Also, let $\lambda_{nn}^{(\alpha,\beta)}$ be the maximum non-trivial zero of the polynomial ${\cal E}_{n0}^{(\alpha, \beta)}(t)$. 

Considering $\beta=\omega\alpha$, we introduce parameter $\omega$ that sets the abscissa of the weight function maximum, $t=\mbox{ln}(1+\omega)$. Then, by subjecting choice of $\alpha_n$ to the requirements
\begin{equation}
\label{A1}
\gamma_n\leq 1\quad \mbox{and} \quad \gamma_n=\underset{\alpha_n}\max\hspace{2pt} \lambda_{nn}^{(\alpha_n, \omega_n\alpha_n)},
\end{equation}
we introduce parameter $\gamma_n$.
In general, equality in (\ref{A1}) can be reached if $\alpha_n,\beta_n \in {\mathbb R}$; inequality holds if $\alpha_n,\beta_n$ are selected from  the set of whole numbers or from a set of rational numbers of interest. Accordingly, for $\gamma_n=1$ we define functions
\[  
\mbox{\boldmath ${\cal Z}$}_{n}^{(\omega_n)}(t)=\{1, {\cal Z}_{nk}^{(\omega_n)}(t)\}_{k=1}^n,\quad
{\cal Z}_{nk}^{(\omega_n)}(t)={\cal E}_{nk}^{(\alpha_n, \omega_n\alpha_n)}(t),
\]
and, for $\gamma_{n}<1$,
\[  
\mbox{\boldmath $\widetilde{{\cal Z}}$}_{n}^{(\omega_n)}(t)=\{1, \widetilde{{\cal Z}}_{nk}^{(\omega_n)}(t)\}_{k=1}^n,\quad
\widetilde{{\cal Z}}_{nk}^{(\omega_n)}(t)={\cal E}_{nk}^{(\alpha_n, \omega_n\alpha_n)}(\gamma_n t).
\] 
Similarly, we define sequences ${\cal Z}_{n0}^{(\omega_n)}(t)$ and $\widetilde{{\cal Z}}_{n0}^{(\omega_n)}(t)$. Choice of $\omega_n$ sets distributions of the zeros of the functions under consideration, 
and ${\cal Z}_{n0}^{(\omega_n)}(1)=\widetilde{{\cal Z}}_{n0}^{(\omega_n)}(1)=0$. System of functions $\left(\mbox{\boldmath $\widetilde{{\cal Z}}$}_{n}^{(\omega_n)} \cup \widetilde{{\cal Z}}_{n0}^{(\omega_n)}\right)(t)$
is easier to construct, in particular for $\alpha_n,\beta_n \in {\mathbb N}$.

Solving the initial value problem, we follow the approach analogous to (\ref{eq5}),(\ref{eq6}). Precisely, we approximate the right-hand side $f(t,y(t)) \equiv g(t)$ of the problem by the functions 
$\mbox{\boldmath $\widetilde{{\cal Z}}$}_{n}^{(\omega_n)}(t)$ making use of calculation of
$g(t)$ on the interval $[0,1]$ at the zeros of $\widetilde{{\cal Z}}_{n0}^{(\omega_n)}(t)$ and at $t=0$. 

For $\alpha_n,\beta_n \in {\mathbb N}$, $\omega_n=1$, and $n=2$ we obtain 
\[
\alpha_2=\beta_2=6,\quad \gamma_2=\ln((15+\sqrt{15})/7),
\]
\[
g(t)\approx g(0)e_{20}(t)+
g(\nu_{1})e_{21}(t)+g(1)e_{22}(t), 
\]
\[
\nu_{1}=1-\mu_{1}/\gamma_2, \quad \mu_{1}=\ln((8+\sqrt{15})/7), 
\]
and the interpolating functions are
\[
e_{20}(t)=1-\frac{30}{7}e^{-\gamma_2 t}+\frac{30}{7}e^{-2\gamma_2 t},
\]
\[
e_{21}(t)=-\sqrt{15}+\frac{15+22\sqrt{15}}{7}e^{-\gamma_2 t}-\frac{15+15\sqrt{15}}{7}e^{-2\gamma_2 t},
\]
\[
e_{22}(t)=+\sqrt{15}+\frac{15-22\sqrt{15}}{7}e^{-\gamma_2 t}-\frac{15-15\sqrt{15}}{7}e^{-2\gamma_2 t}.
\]
Completing construction of the procedure, we find the 
Butcher tableau \footnote[6]{{\sc J.~C. Butcher},
{\em Numerical Methods for Ordinary Differential Equations}, 
John Wiley {\&} Sons, Chichester, 2003.}
\begin{center}
\begin{tabular}{ c|c c c }
  $0    $ & $ $ & $ $ & $ $\\
  $\nu_1$ & $\nu_1-\hat{q}/\gamma_2$ & $-r\nu_1-\hat{s}/\gamma_2$ 
          & $r\nu_1+(\hat{q}+\hat{s})/\gamma_2$\\
  $1    $ & $1-q/\gamma_2$ & $-r+(q+s)/\gamma_2$ 
          & $r-s/\gamma_2$\\ \hline
  $     $ & $1-q/\gamma_2$ & $-r+(q+s)/\gamma_2$ 
          & $r-s/\gamma_2$\\
\end{tabular}
\end{center}
where $q=(8+\sqrt{15})/14$, $\hat{q}=(8-\sqrt{15})/14$, $r=\sqrt{15}$, $s=~3(1+~\sqrt{15})/4$, $\hat{s}=~3(1-~\sqrt{15})/4$.

Further, we derive the stability function 
\[
R_{2}(z)=\frac{1+Az+Bz^2}{1-(1-A)z+Cz^2},
\]
\[
A=1+(3-2r\nu_1)/(2\gamma^2_2),
\] 
\[
B=s\nu_1/(3\gamma_2)-2\hat{s}/(21\gamma_2)+1/(2\gamma_2^2),
\]
\[ 
C=-qr\nu_1/\gamma_2+\hat{q}r/\gamma_2+1/(2\gamma_2^2),
\] 
$\lim_{z\rightarrow -\infty}R_{2}(z)\approx0.543836$, 
and come to the conclusion that the procedure is A-stable. 

Next result represents reconstruction of approximation (\ref{eq15}) that leads to a different type of the stability function. 

The collocation point  $t=0$ is necessary to employ for initiation an explicit algorithm, but for an implicit discretization we can drop direct calculation of the derivative of the solution at the left end of the integration interval. Realizing this approach and making use of ${\cal E}_{nk}(t)$ for constructing the approximation we obtain for $n=~2$ the tableau
\begin{center}
\begin{tabular}{ c| c c }
  $\nu_1$ & $q\nu_1+r/\beta_2$ & $s\nu_1-r/\beta_2$\\
  $1$         & $q+r/\beta_2         $ & $s-r/\beta_2         $\\ \hline
  $ $         & $q+r/\beta_2         $ & $s-r/\beta_2         $\\
\end{tabular}
\end{center}
where
$\beta_2=\ln(3+\sqrt{3})$,
$\nu_1=1-\mu_1/\beta_2$, 
$\mu_1=\ln(2+\sqrt{3})$,  
$q=(3-\sqrt{3})/6$,
$r=\sqrt{3}/6$,
$s=q+2r$.
 
We also find that
\[
R_{2}(z)=\frac{1+Az}{1-(1-A)z+Bz^2},
\]
$A=q\mu_1/\beta_2$, $B=r\mu_1/\beta_2{^2}$, and the procedure is L-stable.
}

One may mention that discretization (\ref{eq15}) contains quadratures with weights of mixed signs, while the procedures outlined in this appendix hold the quadratures with positive weights.


\begin{thebibliography}{1}
{\small
\bibitem{VC1}
{\sc V.~S. Chelyshkov},
Sequences of exponential polynomials, which are orthogonal on the semi-axis, 
{\em Transactions of the Academy of Sciences of the UkSSR, (Doklady AN UkSSR)}, ser. A, No.1, 1987, pp.~14--17. 
(In Russian)
\bibitem{VC2}
{\sc V.~S. Chelyshkov, V.~T. Grinchenko, C.~Liu}, 
Local approach for modeling of quasi-regular 
flows, 
{\em Proceedings of The First AFOSR International Conference on Direct and Large Eddy Simulation,} C. Liu, Z. Liu (editors), 
Greyden Press, Columbus, 1997, pp.~265--272.
\bibitem{Gear}
{\sc C.~W. Gear}, 
{\em Numerical Initial Value Problem in Ordinary Differential Equations}, Englewood Cliffs, NJ: Prentice Hall, 1971. 
\bibitem{VC3}
{\sc V.~T. Grinchenko, V.~S. Chelyshkov},
Transition in the case of low free stream turbulence,
{\em Progress and Challenges in CFD, Methods and Algorithms}, AGARD-CP-578, pp. 26-1--26-9.
\bibitem{Lanczos}
{\sc C. Lanczos}, 
{\em Applied Analysis}, Prentice Hall, Inc., 1956.
\bibitem{Press}
{\sc W.~H. Press, S.~A. Teukolsky, V.~T. Wetterling, and 
B.~P. Flannery}, 
{\em Numerical Recipes}, Cambridge University Press, 1992.
\bibitem{Stoer}
{\sc J. Stoer and R. Bulirsch}, 
{\em Introduction to Numerical Analysis}, 
Springer-Verlag, 1980.
}
\end{thebibliography}
\end{document}